\documentclass[12pt]{article}

\usepackage[english]{babel}
\usepackage{upref,amsfonts,amsxtra,latexsym,color}
\usepackage{amssymb,amsmath,amsthm,epsf,a4wide,mathrsfs,verbatim,hyperref}
\usepackage[all]{xy}

\newtheorem{theorem}{Theorem}[subsection]

\newtheorem{corollary}[theorem]{Corollary}
\newtheorem{proposition}[theorem]{Proposition}

\theoremstyle{definition}
\newtheorem{definition}[theorem]{Definition}
\newtheorem{remark}[theorem]{Remark}

\numberwithin{equation}{section}
\numberwithin{theorem}{section}

\newcommand{\N}{{\mathbb N}}

\newcommand{\R}{{\mathbb R}}
\newcommand{\Cs}{{$C^*$-al\-ge\-bra}}

\newcommand{\Aff}{\mathrm{Aff}}

\date{}

\title{The Space of Tracial States on a $C^*$-Algebra}

\author{Bruce Blackadar and Mikael R\o rdam\footnote{The second named author was supported by a research grant from the Danish Council for Independent Research, Natural Sciences} }

\begin{document}

\maketitle

\begin{center}
\emph{Dedicated to the memory of Bent Fuglede}
\end{center}

\begin{abstract} 
\noindent 
We give a simple and elementary proof that the tracial state
space of a unital \Cs{} is a Choquet simplex, using the center-valued trace on a finite von Neumann algebra.
\end{abstract}

\section{Introduction}  

The object of this paper is to give a simple and elementary proof (modulo well-known facts about
von Neumann algebras) of the following theorem:

\begin{theorem}\label{thm:csimplex}
Let $A$ be a unital \Cs.  Then the tracial state space $T(A)$ of $A$ is a (possibly empty) Choquet simplex.
\end{theorem}

\noindent
This theorem was first proved (in disguised form) in \cite{Thoma:Gruppen}, and another (slightly different)
proof appeared in \cite[Theorem 3.1.18]{Sak:C*-W*}, using the Radon-Nikodym Theorem for operators, or see \cite{Wang}
for a more detailed exposition.  See \cite[p.~120]{Takesaki} for a more ge\-ne\-ral result; the proof
there, specialized to finite von Neumann algebras, gives a proof of Theorem \ref{thm:csimplex}
(essentially the proof from \cite{Thoma:Gruppen}).  There is another somewhat more elementary 
approach due to Pedersen (\cite{Pedersen}, cf.\ \cite[2.8]{CuntzP}).  
There is also a partial converse: every metrizable 
Choquet simplex occurs as the trace space of a separable \Cs, even a simple (unital)
AF algebra.  See also \cite{IoanaSV} and \cite{OrovitzSV}, where it is shown that the trace simplex of certain universal free product \Cs s is the Poulsen simplex. 

The origin of this paper was the experience of both authors that while Theorem \ref{thm:csimplex}
is well known and considered important by operator algebraists, the existing proofs are
not widely known or understood (the result has almost reached the status of folklore). We hope with this exposition to
make the proof of this important result more accessible. Theorem \ref{thm:csimplex} plays a particularly important role in the classification program for simple \Cs s.

Our proof is very simple and elementary given well-known facts about Choquet simplexes and
von Neumann algebras, notably the center-valued trace on finite von Neumann algebras.
The existence of the center-valued trace is a nontrivial but well-known and fundamental result
in the theory of von Neumann algebras, and is in our opinion a reasonable starting point for the
proof of Theorem \ref{thm:csimplex}.  It should be noted that the previous proofs are
somewhat more elementary from first principles since they do not (explicitly) use the center-valued trace.
Deep down all the proofs (except perhaps the Pedersen argument) are really essentially the same, reducing the statement to the fact that
the state space of any commutative von Neumann algebra, indeed any commutative unital 
\Cs, is a Choquet simplex (even a Bauer simplex); the ingredients are largely the same,
but our argument buries all the technicalities in the center-valued trace.
A good exposition of the center-valued trace can be found in \cite{KadRin:opII}.

In Section~\ref{sec:Choquet} we review the basic theory of Choquet simplexes, including a selection of the many characterizations of these.  Using one of these characterizations, we show that a quotient of  a Choquet simplex with an affine cross section is again a Choquet simplex. We use this result in  Section~\ref{sec:main} to prove our main result.  We emphasize that most of the material in Section~\ref{sec:Choquet} is not needed for the proof of Theorem~\ref{thm:csimplex}, and it is included mostly for expository reasons to remind the reader of some of the rich and beautiful theory behind Choquet simplexes. 

We thank George Elliott for helpful comments, and for the reference to \cite{Wang}, and Leonel
Robert for the references to \cite{Pedersen} and \cite{Robert}.

\section{Choquet Simplexes} \label{sec:Choquet}

As mentioned already in the introduction, we review here some of the many (well-known) equivalent conditions for a compact convex subset $K$ of a locally convex topological space $E$ to be a Choquet simplex, stated formally as Theorem~\ref{thm:simplex} below.   Other equivalent conditions are known,
such as the one in \cite{NamiokaP}; this paper also discusses tensor products of Choquet
simplexes, with some interesting questions with applications in quantum information theory.
Incidentally, our Proposition \ref{prop:lift} gives a very simple proof of the equivalence of {\em simplex}
and {\em simplex-like} in that paper.

The \emph{barycenter}  of a probability measure $\mu$ on $K$ is $\int_K x \, d\mu(x)$, which is a point in $K$. (All measures here are understood to be Borel measures.) A probability measure $\mu$ on $K$ is said to be a \emph{boundary measure} if it is a maximal in the sense that if $\nu$ is another probability measure on $K$ satisyfing $\mu(f) \le \nu(f)$, for all $f \in P(K)$, the cone of all continuous convex functions $f \colon K \to \R$, then $\nu = \mu$, cf.\ \cite[Proposition I.4.5]{Alf:convex}. If $K$ is metrizable, then the set $\partial_e K$ of extreme points of $K$ is a $G_\delta$-set, and $\mu$ is a boundary measure if and only if $\mu(\partial_e K)=1$. Each $x \in K$ is the barycenter of some boundary measure on $K$ (Choquet--Bishop--de Leeuw, cf.\ \cite[Theorem I.4.8]{Alf:convex}). 

Let $\Aff(K)$ denote the set of continuous affine functions $f \colon K \to \R$, and equip $\Aff(K)$ with the standard pointwise ordering: $f \ge 0$ iff $f(x) \ge 0$ for all $x \in K$.  ($\Aff(K)$ also has a
{\em strict ordering}, where $f\gg 0$ iff $f=0$ or $f(x)>0$ for all $x\in K$, which is important in
applications to the classification program.)

The inclusion $K \subseteq E$ is said to be \emph{regular} if the affine subspace of $E$ spanned by $K$ is a hyperplane in $E$ not containing $0$, in other words, equal to $\varphi^{-1}(\{1\})$, for some $\varphi \in E^*$. Let $\tilde{K}$ denote the cone $\{\alpha x : x \in K, \alpha \ge 0\}$ in $E$, and let $(E, \tilde{K})$ denote $E$ equipped with the order relation given by $x \le y$ iff $y-x \in \tilde{K}$, for $x,y \in E$.   Every compact convex set $K$ is isomorphic (affinely homeomorphic) to an essentially unique regularly embedded compact convex set: the subset of $E=\Aff(K)^*$ consisting of the states
(order-preserving linear maps $\phi \colon \Aff(K)\to\R$ with $\phi(u)=1$, where $u$ is the constant function 
1) is a regularly embedded compact convex set isomorphic to $K$.

We first state a proposition which is a combination of results in \cite{Alf:convex} and \cite{Goo:partial}.

\begin{proposition} \label{prop:Riesz}
Consider the following properties of a partially ordered real vector space $E$:
\begin{enumerate}
\item $E$ is a lattice.
\item $E$ has the \emph{Riesz interpolation property}, i.e., for all $x_1,x_2, y_1, y_2 \in E$ with $x_j \le y_i$, for $i,j=1,2$, there exists $z \in E$ with $x_j \le z \le y_i$, for $i,j=1,2$.
\item $E$ has the \emph{Riesz decomposition property}, i.e., whenever $0 \le y \le x_1+x_2$ in $E$, then there exist $0 \le y_j \le x_j$, $j=1,2$, in $E$ such that $y = y_1+y_2$.
\item Whenever $x_1, \dots, x_n, y_1, \dots, y_m$ are positive elements in $E$ with $\sum_{i=1}^n x_i = \sum_{j=1}^m y_j$, then there exist $z_{ij} \ge 0$ in $E$ such that $x_i = \sum_{j=1}^m z_{ij}$ and $y_j = \sum_{i=1}^n z_{ij}$, for all $1 \le i \le n$ and $1 \le j \le m$.
\end{enumerate}
Then {\rm{(i)}} $\Rightarrow$ {\rm{(ii)}} $\Leftrightarrow$ {\rm{(iii)}} $\Leftrightarrow$ {\rm{(iv)}}, and {\rm{(ii)}} $\Rightarrow$ {\rm{(i)}} if $E^+$, the positive cone of $E$, is locally compact in some locally convex topology on $E$. In particular, if $K$ is a compact convex set regularly embedded in $E$ and $E^+ = \tilde{K}$ is the positive cone in $E$ generated by $K$, as above, then all four conditions above are equivalent.
\end{proposition}

\medskip \noindent It is clear that (i) $\Rightarrow$ (ii) as we may take $z = x_1 \vee x_2$. 
For proofs of the equivalence of (ii) and (iii) and the implication (ii) $\Rightarrow$ (i) (assuming locally compactness of $E^+$), see, e.g., Alfsen, \cite[Proposition II.3.1 and Proposition II.3.2]{Alf:convex}. If $K$ is compact and convex, then $\tilde{K}$ is locally compact, since each point in $\tilde{K}$ belongs to the interior of the compact set $\{tx : 0 \le t \le r, x \in K\}$, for some $r >0$. 

For the convenience of the reader, and because this point is central for the proof of Proposition~\ref{prop:lift} below, and hence for our arguments in the next section, we review here the proof of (iii) $\Leftrightarrow$ (iv). 

(iv) $\Rightarrow$ (iii). Assume  $0 \le y \le x_1+x_2$. Set $y_1=y$ and $y_2=x_1+x_2-y$. Then $x_1+x_2=y_1+y_2$, and there exist $z_{ij} \ge 0$ in $E$ such that $x_i=z_{1i}+z_{2i}$ and $y_j = z_{j1}+z_{j2}$, for $i,j =1,2$. In particular, $y = z_{11}+z_{12}$ and $0 \le z_{1i} \le x_i$, as wanted.

(iii) $\Rightarrow$ (iv). We first show, by induction, that whenever  $n \ge 2$ and $0 \le y \le \sum_{j=1}^n x_j$ in $E$, there exist $0 \le y_j \le x_j$, $1 \le j \le n$, in $E$ such that $y=\sum_{j=1}^n y_j$. 
Indeed, the base step  $n=2$ holds by assumption. Let $n > 2$ and suppose the claim has been verified for $n-1$. 
 Set $x_1' = \sum_{j=1}^{n-1} x_j$ and $x_2' =x_n$. Then $y \le x_1'+x_2'$ and so there exist $0 \le y_j' \le x_j'$, $j=1,2$,  satisfying $y=y_1'+y_2'$. Set $y'=y_1'$ and $y_n =y_2'$. Then $0 \le y' \le  \sum_{j=1}^{n-1} x_j$. By the induction hypothesis there exist $0 \le y_j \le x_j$, $1 \le j \le n-1$, in $E$ such that $y'=\sum_{j=1}^{n-1} y_j$, whence $y = \sum_{j=1}^n y_j$. This completes the induction step.

We proceed to prove (iv) by induction on $m \ge 1$. The base step $m=1$ is trivial. Let $m \ge 2$ and assume the claim holds for $m-1$. Let  $x_1, \dots, x_n \ge 0$ and $y_1, \dots, y_m \ge 0$ satisfy $\sum_{i=1}^n x_i = \sum_{j=1}^m y_j$. Then $0 \le y_m \le \sum_{i=1}^n x_i$, so there exist $0 \le z_{i,m} \le x_i$, $1 \le i \le n$, in $E$ such that $y_m = \sum_{i=1}^n z_{i,m}$. Set $x_i' = x_i - z_{m,i} \ge 0$. Then $\sum_{i=1}^n x_i' = \sum_{j=1}^{m-1} y_j$. By the induction hypothesis we can find $z_{ij} \ge 0$, $1 \le i \le n, 1 \le j \le m-1$,  in $E$ such that $$x_i' = \sum_{j=1}^{m-1} z_{ij}, \qquad y_j = \sum_{i=1}^n z_{ij}, \quad 1 \le i \le n, \; 1 \le j \le m-1.$$ The former identity implies  $x_i = x_i' + z_{i,m} = \sum_{j=1}^{m} z_{ij}$, as desired. \hfill $\square$
 
\bigskip  \noindent We shall now present a selection of equivalent conditions characterizing Choquet simplexes. These are all classical results which can be found in several textbooks (see, e.g.,  the comments below the theorem), and at least the equivalence of (i) and (ii) goes back to Choquet.  See also the notes at the end of Chapter II.3 of \cite{Alf:convex} for a more detailed account on the history. 

A \emph{homothetic image} of $K$ in $E$ is a set of the form $\alpha K + x$, where $\alpha  >0$ and $x \in E$.

Finally, for ease of notation, when we say that $\sum_{j=1}^n \alpha_j x_j$ is a \emph{proper convex combination} in $K$, we mean that $\alpha_1, \dots, \alpha_n > 0$, $x_1, \dots, x_n \in K$ and $\sum_{j=1}^n \alpha_j=1$. 

\begin{theorem}
\label{thm:simplex} Let $K$ be a compact convex set. The following conditions are equivalent. 
\begin{enumerate}
\item Each point in $K$ is the barycenter of a unique boundary measure  on $K$. 
\item The ordered vector space $\Aff(K)^*$ is a lattice.
\item {\emph{(Riesz Decomposition Property)}} Whenever $x \in K$ in two different ways is a proper convex combination in $K$, i.e.,
$$x = \sum_{j=1}^n \alpha_j x_j = \sum_{i=1}^m \beta_i y_i,$$
then  there exists a new convex combination $x = \sum_{j=1}^n \sum_{i=1}^m \gamma_{ij} z_{ij}$ in $K$, satisfying
\begin{equation} \label{star}
\alpha_j = \sum_{i=1}^m \gamma_{ij}, \quad \beta_i = \sum_{j=1}^n \gamma_{ij}, \qquad
x_j = \sum_{i=1}^m \alpha_j^{-1} \gamma_{ij} z_{ij}, \quad y_i = \sum_{j=1}^n \beta_i^{-1} \gamma_{ij} z_{ij},
\end{equation} 
for all $1 \le j \le n$ and $1 \le i \le m$ 
\item The ordered vector space $\Aff(K)$ has the Riesz Interpolation Property in either (both) the
ordinary and strict orderings (i.e., if $f_1,f_2\leq g_1,g_2$, then there is an $h$ with
$f_1,f_2\leq h\leq g_1,g_2$, and similarly for $\ll$).
\item Whenever $K$ is embedded in a locally convex vector space $E$, the intersection of $K$ with any homothetic image of itself is either empty, a point, or another homothetic image of $K$.
\end{enumerate}
\end{theorem}

\noindent Condition (ii) can further be reformulated as described in Proposition~\ref{prop:Riesz}.

Proofs of the equivalences of (i), (ii) and (iii) can be found, e.g., in Alfsen, \cite[Propositions II.3.3 and II.3.6]{Alf:convex}.  The equivalence of (ii) and (iv) is proved in \cite[Theorem 11.4]{Goo:partial}; see also \cite[Proposition 10.7]{Goo:partial} for a proof of (ii) $\Rightarrow$ (iii). 
The beautiful geometric interpretation of simplexes in (v) is portrayed on the cover of Phelps' book, \cite{Phelps:Choquet}; it is stated but not proved in \cite{Phelps:Choquet} that (v) is equivalent to being a Choquet simplex. One can fairly easily see that if (v) holds in some representation of $K$ in a locally convex vector space $E$, then it holds for all such representations. 

For completeness of exposition we include here the short proof of (ii) $\Leftrightarrow$ (iii) via Proposition~\ref{prop:Riesz}, and, for the record, also the (not so short) proof of (ii) $\Leftrightarrow$ (v). We shall use the following setup. Retaining the notation given above Proposition~\ref{prop:Riesz}~(iv), let $K \subseteq E$ be a regular inclusion, and $\tilde{K} \subseteq E$ the positive cone in $E$ generated by $K$. A canonical choice of such an inclusion is obtained by taking $E =\Aff(K)^*$. Condition (ii) then says that $(E, \tilde{K})$ is a lattice. Fix $\varphi \in E^*$ such that $K \subseteq \varphi^{-1}(1)$. Observe that if $e \in E$ is positive, then $\varphi(e) \ge 0$; $\varphi(e) = 0$ if and only if $e=0$; and $e \in K$ if and only if $\varphi(e)=1$.  

Proposition~\ref{prop:Riesz}~(iv) $\Rightarrow$ (iii):  Let $x = \sum_{j=1}^n \alpha_j x_j = \sum_{i=1}^m \beta_i y_i$ be two proper convex combinations in $K$. Set $\bar{x}_j = \alpha_j x_j$ and $\bar{y}_i = \beta_i y_i$. Then there exist positive elements $\bar{z}_{ij} \in E$ satisfying $\bar{x}_i = \sum_{j=1}^m \bar{z}_{ij}$ and $\bar{y}_j = \sum_{i=1}^n \bar{z}_{ij}$. Set $\gamma_{ij} = \varphi(\bar{z}_{ij}) \ge 0$.  Then $\sum_{j=1}^m \gamma_{ij} = \varphi(\bar{x}_i) = \alpha_i$ and $\sum_{i=1}^n \gamma_{ij} = \varphi(\bar{y}_j) = \beta_j$. Hence \eqref{star} holds with $z_{ij} = \gamma_{ij}^{-1} \bar{z}_{ij} \in K$, if $\gamma_{ij} \ne 0$, and any $z_{ij} \in K$ when $\gamma_{ij} = 0$. 

(iii) $\Rightarrow$ Proposition~\ref{prop:Riesz}~(iv): Let $\bar{x}_1, \dots, \bar{x}_n, \bar{y}_1, \dots, \bar{y}_m$ be positive elements in $E$ satisfying $x:= \sum_{i=1}^n \bar{x}_i = \sum_{j=1}^m \bar{y}_j$.
Upon multiplying these equations by a suitable positive number, we may assume that $\varphi(x) = 1$, i.e., $x \in K$. We may also assume that all $\bar{x}_i$ and $\bar{y}_j$ are non-zero. Put $\alpha_i = \varphi(\bar{x}_i)$  and $\beta_j = \varphi(\bar{y}_j)$, and set $x_i = \alpha_i^{-1} \bar{x}_i$ and $y_j = \beta_j^{-1} \bar{y}_j$. Then $x_i,y_j \in K$ because $\varphi(x_i)=\varphi(y_j)=1$, and  $x = \sum_{i=1}^n \alpha_i x_i = \sum_{j=1}^m \beta_j y_j$ are proper convex combinations in $K$. Accordingly, there exist $z_{ij} \in K$ and $\gamma_{ij} \ge 0$ such that \eqref{star} holds. Set $\bar{z}_{ij} = \gamma_{ij} z_{ij} \ge 0$. Then $\bar{x}_i = \sum_{j=1}^m \bar{z}_{ij}$ and $\bar{y}_j= \sum_{i=1}^n \bar{z}_{ij}$, as desired.

For the proof of (ii) $\Leftrightarrow$ (v), note first that for $a, b, c \in E$,
$$a \le b \iff b \in \widetilde{K}+a, \qquad c = a \vee b \iff (\widetilde{K} + a) \cap (\widetilde{K}+b) = \widetilde{K}+c.$$ We shall make frequent use of the following identities:
\begin{equation} \label{eq:2}
 \widetilde{K} \cap \varphi^{-1}(r) = rK,
\end{equation}
for all $r \ge 0$, and 
\begin{equation} \label{eq:3}
(\widetilde{K} +a) \cap \varphi^{-1}(r) = (r-\varphi(a))K +a
\end{equation}
for  $a \in E$ and  $r \ge \varphi(a)$.

(ii) $\Rightarrow$ (v): Let $r > 0$ and $a \in E$ be given, and suppose that  $K \cap (rK+a)$  is non-empty. Then $1=r+\varphi(a)$. Set $c = 0 \vee a$. Then $\widetilde{K} \cap (\widetilde{K}+b) = \widetilde{K} + c$. 
By \eqref{eq:2} and \eqref{eq:3} we therefore get
$$K \cap (rK+b) = \widetilde{K} \cap (\widetilde{K} +b) \cap \varphi^{-1}(1) = (\widetilde{K} +c) \cap \varphi^{-1}(1) = tK + c,$$
with $t = 1 - \varphi(c)$. (If $\varphi(c) > 1$, then $K \cap (rK+b)  = \emptyset$.)

(v) $\Rightarrow$ (ii). We asssume in the proof that $K$ is not a singleton. If it were, then the ordered vector space $E$ would be $\R$ with the usual ordering, which is a lattice. Secondly, if (v) holds, then it also holds that the intersection of any two homothetic images of $K$, if not empty or a point, is another homothetic image of $K$.

Let $a,b \in E$. For each $t \in \R$, let $C_t$ be the set of $c \in E$ for which $a,b \le c$ and $\varphi(c) = t$. Note that $C_t$ is a convex set, and that
\begin{equation} \label{eq:C}
C_t = (\tilde{K}+a) \cap (\tilde{K}+b) \cap \varphi^{-1}(t) = \big((t-\varphi(a))K + a\big) \cap 
\big((t-\varphi(b))K+b\big),
\end{equation}
when $t \ge \varphi(a), \varphi(b)$. Set $I= \{t \in \R : C_t \ne \emptyset\}$. By assumption, and since $C_t \subseteq \varphi^{-1}(t)$, there exists $c_t$ (a priori depending on $t$) such that 
\begin{equation} \label{eq:B}
C_t = (t-\varphi(c_t))K + c_t, \qquad t \in I.
\end{equation}
 
We proceed to prove the following three claims, which will complete the proof.
\begin{itemize}
\item[(a)] $C_t + sK \subseteq C_{t+s}$, for all $s \ge 0$ and $t \in I$. 
\item[(b)] If $tK+e \subseteq t'K + e'$, for some $t,t' \ge 0$ and $e,e' \in E$, then $t \le t'$ and $e' \le e$; and if in addition $t=t'$, then $e=e'$.
\item[(c)] $c_t = c$ is independent of $t \in I$, and $c = a \vee b$. 
\end{itemize}

(a) is trivial: if $x \in C_t + sK$, then $x \ge a,b$ and $\varphi(x) = t+s$, so $x \in C_{t+s}$.

(b).  Let $e,e' \in E$ be such that $tK+e \subseteq t'K + e'$, for some $t,t' \ge 0$.  We first show that $t 
\le t'$. Suppose otherwise. Then, for all $m \ge 1$,
$$t' K + me' \supseteq (t-t')K + t'K +e+(m-1)e'  \supseteq \cdots \supseteq m(t-t')K+t'K +me.$$
Pick $x \in K$ such that $y:= (t-t')x + e-e' \ne 0$. Then $\{my +t'x : m \in \N\} \subseteq t'K$, 
which is impossible, as $t' K$ is compact.  

Suppose next that $t=t'$. Then $tK + (e-e') \subseteq tK$, which entails $tK + m(e-e') \subseteq tK$, for all integers $m \ge 1$. As $tK$ is compact, this implies $e=e'$. 

Consider  finally the case $t < t'$. Set $r' = \varphi(e'), r = \varphi(e)$, and $s = t+r=t'+r'$. Note that $r' < r \le s$. Assume, to reach a contradiction, that $e' \nleq e$. Then  $e \notin   (r-r')K +e'$.  Choose $y$ in the compact convex set $(r-r')K+e'$ so that  
$$\{\alpha y +(1-\alpha)e : 0 \le \alpha < 1\} \cap \big((r-r')K+e'\big) = \emptyset.$$
(One can find such a $y$ in the line segment joining any point in $K$ to $e$.) Write $y = (r-r')x+e'$, with $x \in K$, set
$$z = tx +e \in tK+e \subseteq t'K + e' \subseteq \varphi^{-1}(s), \qquad \alpha = \frac{s-r}{s-r'}, \quad 1-\alpha = \frac{r-r'}{s-r'},$$
 and set $u = \alpha e' +(1-\alpha) z$. Then $ u \in (r-r')K + e'$ because $z \in  t'K + e'$ and $(1-
 \alpha) t' = r-r'$. On the other hand,
 $$
 u= \alpha e' + (1-\alpha)tx + (1-\alpha)e = \alpha \big(e' + (r-r')x\big) + (1-\alpha)e = \alpha y + (1-\alpha)e,
 $$
 which contradicts the choice of $y$. This completes the proof of (b).

(c). It follows from Equations~\eqref{eq:C} and \eqref{eq:B} and (b) that  $a,b \le c_t$, for all $t \in I$. 
If $t'\le t$ are in $I$, then by (a),
\begin{eqnarray} \label{eq:A}
(t-\varphi(c_{t'}))K + c_{t'} = C_{t'} + (t-t')K \subseteq C_{t} = (t-\varphi(c_{t}))K + c_{t},
\end{eqnarray}
so $c_{t} \le c_{t'}$ by (b). Set $t_0 = \inf I \ge \max\{\varphi(a), \varphi(b)\}$. We claim that $\varphi(c_t) = t_0$, for all $t \in I$. Since $c_t \ge a,b$, it follows that $\varphi(c_t) \in I$, so $\varphi(c_t) \ge t_0$. 
Suppose $t \in I$ and $t_0 < \varphi(c_t)$. Take $t' \in I$ with $t_0 < t' < \varphi(c_t) \le t$. Then $c_t \le c_{t'}$, so $\varphi(c_t) \le \varphi(c_{t'}) \le t'$, a contraction.  By \eqref{eq:A} and (b) this further implies that $c_t$ are equal, hence equal to $c$, for some $c \in E$.  We already noted that $c \ge a,b$. Suppose $d \in E$ also satisfies $d \ge a,b$. Then $d \in C_{\varphi(d)}$, which implies $d \ge c_{\varphi(d)} = c$. This proves that $c = a \vee b$. 
\hfill $\square$

\begin{definition} A compact convex set $K$ satisfying the equivalent conditions of Theorem~\ref{thm:simplex} is a {\em Choquet simplex}.
\end{definition}

\noindent
An $n$-simplex in Euclidean space is a Choquet simplex with $n+1$ extreme points.
For the purposes of Theorem \ref{thm:csimplex}, we will regard $\emptyset$ (which technically
satisfies the conditions of Theorem~\ref{thm:simplex}) as a Choquet simplex (a (-1)-simplex?)

The set of extreme points of a Choquet simplex need not be closed, although the set of extreme
points of a metrizable Choquet simplex is a $G_\delta$-set and thus a Polish space (the set of extreme
points of a nonmetrizable Choquet simplex need not even be a Borel set).  A Choquet simplex
in which the set of extreme points is closed is called a {\em Bauer simplex}.

The Riesz Interpolation Property (Theorem~\ref{thm:simplex}(iv)) is a weak lattice property.  But if $K$ is a 
Choquet simplex, $\Aff(K)$ is not a lattice in general; in fact, $\Aff(K)$ is a lattice if and only if
$K$ is a Bauer simplex \cite[Theorem II.4.1]{Alf:convex}.

The plural of {\em simplex} is often written {\em simplices}, which is correct Latin, but {\em simplexes}
seems to be more common now, and is more in accord with the plural of {\em complex}
(as in {\em CW-complex}).

\medskip
\begin{remark} \label{rem:Riesz} Observe that Theorem~\ref{thm:simplex}(iii) is a strengthening of the statement that each element in $\mathrm{conv}(\partial_e K)$ has a unique representation as a convex combination of extreme points in $K$. Indeed, if $x_1, \dots, x_n$, respectively, $y_1, \dots, y_m$, are distinct extreme points of $K$ and $x = \sum_{i=1}^n \alpha_i x_i = \sum_{j=1}^m \beta_j y_j$ are convex combinations in $K$. Then \eqref{star} holds for some $z_{ij} \in K$ and $\gamma_{ij} \ge 0$. If $\gamma_{ij} \ne 0$, then $x_i = z_{ij} = y_j$. This implies $m=n$ and that $\{x_1, \dots, x_n \} = \{y_1, \dots, y_m\}$. 

A compact convex set in which each element in $\mathrm{conv}(\partial_e K)$ has a unique representation as a convex combination of extreme points in $K$ is not necessarily a Choquet simplex: a counterexample is easy to construct, e.g., dividing out the Bauer simplex $\Delta$ with extreme
points $[0,1]$ (sitting regularly in $\Aff(\Delta)^*$) by the one-dimensional subspace spanned by $\mu-\nu$, where $\mu$ and $\nu$
are distinct continuous probability measures on $[0,1]$. 
\end{remark}

\bigskip 

\noindent Using the characterization of Choquet simplexes in Theorem \ref{thm:simplex}(iii) we can easily prove the following:

\begin{proposition} \label{prop:lift} Let $K$ be a compact convex set and let $\Delta$ be a Choquet simplex.  Assume there exists an affine surjection $\varphi \colon \Delta \to K$ with an affine cross section $\lambda \colon K \to \Delta$, i.e.\ $\varphi\circ\lambda=id_K$. Then $K$ is a Choquet simplex.
\end{proposition}

\noindent Note that we do not require that the maps $\varphi$ and $\lambda$ above are continuous!
(But they must be affine.)

\begin{proof} Let $x \in K$ be written as two proper convex combinations $x = \sum_{j=1}^n \alpha_j x_j = \sum_{i=1}^m \beta_i y_i$ in $K$. Then
$$\lambda(x) = \sum_{j=1}^n \alpha_j \lambda(x_j) = \sum_{i=1}^m \beta_i \lambda(y_i)$$
are proper convex combinations in $\Delta$. 
Since $\Delta$ is a Choquet simplex, there exists a convex combination $\lambda(x) = \sum_{j=1}^n \sum_{i=1}^m \gamma_{ij} \bar{z}_{ij}$ in $\Delta$, satisfying \eqref{star} with $\bar{z}_{ij}$ in the place of $z_{ij}$.
Set $z_{ij} = \varphi(\bar{z}_{ij}) \in K$. Then $x$ is the convex combination $\sum_{j=1}^n \sum_{i=1}^m \gamma_{ij} z_{ij}$ in $K$ and \eqref{star} holds.
\end{proof}

\section{The Main Result} \label{sec:main}

\noindent Let $A$ be a unital \Cs. The set $S(A)$  of all states on $A$ is a compact convex set in the weak$^*$ topology, and $T(A)$, the set of all tracial states on $A$, is a closed convex subset of $S(A)$ (possibly empty). 

We now give the proof of Theorem \ref{thm:csimplex}.

\begin{proof} Consider the inclusion $\iota \colon A \to A^{**}$ of $A$ into its bidual, and recall that $A^{**}$ can be identified with $\Phi_u(A)''$, the von Neumann algebra envelope of $A$ in its universal representation $\Phi_u$, see, e.g., \cite[Section 10.1]{KadRin:opII}. A crucial property of this inclusion is the identification of the dual $A^*$ of $A$ with the predual $(A^{**})_*$ of $A^{**}$. In particular, we have an affine bijection $\rho \mapsto \rho \circ \iota$ from the set of normal states on $A^{**}$ to $S(A)$, which restricts to an affine bijection from the set of normal tracial states on $A^{**}$ onto $T(A)$. 

Write $A^{**} = A^{**}_{\mathrm{fin}} \oplus A^{**}_{\mathrm{inf}}$ as a direct sum of its finite and its properly infinite parts. Since traces on $A^{**}$ obviously vanish on $A^{**}_{\mathrm{inf}}$, we can view these as being defined on $A^{**}$, or on $M:=A^{**}_{\mathrm{fin}}$, as we please.

Let $\mathrm{tr}_c \colon M \to Z:=Z(M)$ denote the (unique) normal center-valued trace, \cite[Theorem 8.2.8]{KadRin:opII}, where $Z$  denotes the center of $M$. We may extend the center-valued trace to all of $A^{**}$ by letting it be zero on  $A^{**}_{\mathrm{inf}}$. Each tracial state on $M$ factors through the center-valued trace (as a consequence of the Dixmier property, see \cite[Theorem 8.3.6]{KadRin:opII}). In summary, we have the following identifications:
\begin{equation} \label{dagger}
 T(A) \longleftrightarrow \{ \text{normal tracial states on $M$}\} \longleftrightarrow \{ \text{normal states on $Z$}\},
 \end{equation}
where the first identification (from right to left) is given by $\tau \mapsto \tau \circ \iota$, for $\tau$ a normal tracial state on $M$,  and the last (again from right to left) is given by $\rho \mapsto \rho \circ \mathrm{tr}_c$, where $\rho$ is a normal state on $Z$. 

To prove the result we use Proposition~\ref{prop:lift} with $K = T(A)$ and $\Delta = S(Z)$, the Bauer simplex of all states on the commutative von Neumann algebra $Z$, which, via Riesz' Representation Theorem can be identified with the Bauer simplex of probability measures on the spectrum of $Z$.  

Let $\varphi \colon S(Z) \to T(A)$ be given by $\varphi(\rho) = \rho \circ  \mathrm{tr}_c \circ \iota$, and let $\lambda \colon T(A) \to  S(Z)$ be the composition of the two (left-to-right) maps in \eqref{dagger} 
above, i.e.,  $\lambda(\tau) = \rho$, where $\rho$ is the unique normal state on $Z$ for which $\tau =  \rho \circ  \mathrm{tr}_c \circ \iota = \varphi(\rho)$. The maps $\varphi$ and $\lambda$ are clearly affine (note that $\lambda$ is not continuous in general), 
and $\varphi \circ \lambda$ is the identity on $T(A)$, by definition of $\lambda$, so it follows from Proposition~\ref{prop:lift} that $T(A)$ is a Choquet simplex.
\end{proof}

\section{The Tracial State Space of a von Neumann Algebra} \label{sec:center}

Let $M$ be a finite von Neumann algebra and consider, as in the proof of Theorem~\ref{thm:csimplex}, the (unique, normal) center-valued trace $\mathrm{tr}_c \colon M \to Z$, where $Z$ is the center of $M$. Combining \cite[Theorem 8.2.8]{KadRin:opII} (existence of center-valued trace) and  \cite[Theorem 8.3.6]{KadRin:opII} (the Dixmier property) we conclude as in the proof of Theorem~\ref{thm:csimplex} that every (normal) tracial state on $M$ is of the form $\rho \circ \mathrm{tr}_c$ for some (normal) state on $Z$. This shows that the set  $T(M)$ of tracial states on $M$ is affinely homeomorphic to the Bauer simplex of all states on the commutative von Neumann algebra $Z$, and that the convex set of normal tracial states on $M$ is affinely homemorphic to the convex set of normal states on $Z$. In particular, the trace simplex of a von Neumann algebra is always a Bauer simplex (unlike the case for general for separable unital \Cs s).  We summarize:

\begin{theorem}\label{thm:vsimplex}
Let $M$ be a finite von Neumann algebra with center $Z$ and center-valued trace 
$\mathrm{tr}_c \colon M \to Z$.  
Then every tracial state $\tau$ on $M$ is uniquely of the form $\tau=\rho\circ\mathrm{tr}_c$, 
where $\rho$
is a state on $Z$.  The trace $\tau$ is normal on $M$ if and only if the state $\rho$ is normal
on $Z$.
\end{theorem}  

Since the center-valued trace is idempotent (a conditional expectation),
the correspondence between $T(M)$ and $S(Z)$ is nothing but restriction.  So we obtain:

\begin{corollary}\label{cor:vsimplex}
Let $M$ be a finite von Neumann algebra with center $Z$.  Then the restriction map gives an affine bijection and weak* homeomorphism from $T(M)$ onto $S(Z)$, i.e., every state on $Z$ extends uniquely to a tracial state on $M$.  This map sends the normal traces on M onto the normal states of Z.  In particular, $T(M)$ is a Bauer simplex.
\end{corollary}

This statement makes no explicit mention of the center-valued trace.  It is slightly less precise than the
statement of Theorem \ref{thm:vsimplex}, but is clean and captures the essence of what is
needed for our proof of Theorem \ref{thm:csimplex}.  We actually do not need the
full force of this corollary for our proof (only the normal part), but
it makes for a nice exposition of the argument.

Of course, the statement that $T(M)$ is a Bauer simplex applies to every von Neumann
algebra $M$, since $T(M)$ is the same as the tracial state space of the finite part of $M$, although
it will be empty if $M$ is properly infinite.

\begin{remark} Let us return to the convex set of normal traces on $M$, which is identified with the set of normal states on $Z(M)$. The commutative von Neumann algebra $Z(M)$ is isomorphic to $L^\infty(\Omega, \mu)$ for some measure space $(\Omega, \mu)$. The set of normal linear functionals on $Z(M)$ is equal to the predual of $Z(M)$, which we identify with $L^1(\Omega,\mu)$. More specifically, each $f \in L^1(\Omega,\mu)$ produces the normal functional $\rho_f(g) = \int_\Omega fg \, d\mu$, for $g \in L^\infty(\Omega,\mu)$.  Let $L^1_\R(\Omega,\mu)$ and $L^1_{\sigma}(\Omega,\mu)$ denote the real-valued functions in $L^1(\Omega,\mu)$, respectively, the positive functions with integral 1, which  correspond to the hermitian normal functionals, respectively, the normal states on $L^\infty(\Omega,\mu)$. The inclusion $L^1_{\sigma}(\Omega,\mu) \subseteq L^1_\R(\Omega,\mu)$ is regular. The cone generated by $L^1_{\sigma}(\Omega,\mu)$ is the set of positive functions, which in turn induces the standard ordering on $L_\R^1(\Omega,\mu)$. It is clear that $L_\R^1(\Omega,\mu)$ is a lattice. 

This shows that the convex set of normal states on $M$ is a Choquet simplex provided that it is compact. More specifically, we would need to equip $L_\R^1(\Omega,\mu)$ with a locally convex topology such that $L^1_{\sigma}(\Omega,\mu)$ becomes compact. This can for example be done whenever $L_\R^1(\Omega,\mu)$  is a dual space, which typically it isn't. This is, however, the case when $M$ is the bidual of a \Cs, and one can use this argument to provide another proof of Theorem~\ref{thm:csimplex}, using (ii) of Theorem~\ref{thm:simplex} rather than (iii), still using the center-valued trace.

The set of normal states of $L^\infty([0,1],\lambda)$ is equal to $L^1_\sigma([0,1],\lambda)$, and it is a well-known fact that this convex set does not have any extreme points, and hence is not compact in any reasonable topology, and certainly not a Choquet simplex.
\end{remark}

\begin{remark}
It would be nice to give a similar proof of the result in \cite{BlackadarH} that the set $QT(A)$ of
normalized \emph{quasitraces} on a unital \Cs~$A$ is a Choquet simplex.  But there is no obvious
analog of the center-valued trace or, indeed, the universal enveloping von Neumann algebra
$A^{**}$.  For any fixed normalized quasitrace $\tau$ on (unital) $A$ there is a finite AW$^*$-algebra
$M$ and an embedding of $A$ as a ``weakly dense'' $C^*$-subalgebra of $M$ (there is no strong
or weak topology in this case, so this must be taken with a grain of salt), such that $\tau$ extends
uniquely to a normal quasitrace on $M$.  $M$ is given by an ultraproduct construction rather
than the GNS construction in the trace case.  One would somehow need to tie all these $M$'s 
together into a single finite AW$^*$-algebra with a ``center-valued quasitrace.''  This might be
possible via some version of the ultraproduct construction.  But the resulting proof would
probably be more complicated than the proof in \cite{BlackadarH}.  (It should be noted that in
the case where $A$ is exact, so $QT(A)=T(A)$, the proof in \cite{BlackadarH} gives an alternate
proof of Theorem \ref{thm:csimplex}.)  See also \cite{Robert} for an adaptation of Pedersen's
approach to proving that $QT(A)$ is a Choquet simplex.
\end{remark}


%
%

\providecommand{\bysame}{\leavevmode\hbox to3em{\hrulefill}\thinspace}
\providecommand{\MR}{\relax\ifhmode\unskip\space\fi MR }
\providecommand{\MRhref}[2]{%
  \href{http://www.ams.org/mathscinet-getitem?mr=#1}{#2}
}
\providecommand{\href}[2]{#2}

\vspace{1cm}

\noindent
Bruce Blackadar \\Department of Mathematics and Statistics \\ University of Nevada, Reno 89557 \\ Email: bruceb@unr.edu \\ https://bruceblackadar.com

\vspace{.6cm}

\noindent
Mikael R\o rdam \\
Department of Mathematical Sciences \\
University of Copenhagen \\ 
Universitetsparken 5, DK-2100, Copenhagen \O \\
Denmark \\
Email: rordam@math.ku.dk\\
WWW: http://web.math.ku.dk/$\sim$rordam/

\end{document}